\documentclass[english]{article}
\usepackage[T1]{fontenc}
\usepackage[latin9]{inputenc}
\usepackage{amsmath}
\usepackage{babel}
\begin{document}
\title{The Local Isometric Embedding of Two-Metrics of Low Differentiability
in Euclidean Three Space.}

\maketitle
Author: Edgar Kann

Affiliation: Queens College, CUNY
\begin{abstract}
We prove that the isometric embedding of any metric of differentiability
class $C^{1}$ in $E^{3}$ exists. We use simplified notation for
the given metric, namely geodesic parameters, and level parameters
for the embedded surface in $E^{3}.$ Central to our discussion will
be solutions of initial value problems for two first order non-linear
partial differential equations. We also make use of the classical
theory of linear algebraic systems. We will prove local isometric
embedding. An example is given for which the Gaussian curvature of
the metric is equal to one but the embedded surface is non-analytic.
\end{abstract}

\section{Introduction.}

The Main Problem is to prove the existence of an isometric embedding
in $E^{3}$ of the given metric
\[
\omega=\bar{E}(\bar{u},\bar{v})d\bar{u}^{2}+\bar{F}d\bar{u}d\bar{v}+\bar{G}(\bar{u},\bar{v})d\bar{v}^{2}
\]

where $\bar{E}=1,$ $\bar{F}$$=0,\bar{G}>0,\bar{G}\epsilon C^{1}.$
We emphasize that this is a local problem. That is, the given metric
is given locally and the existence is to be proved locally.

This is called the $\mathit{geodesic}$ $form\mathit{}$ of the metric.
Any metric can be assumed to be in the geodesic form without loss
of generality \cite{key-3}, \cite{key-4}.

The isometrically embedded surface will be presented in the form of
a vector function 
\[
X(u,v)=x(u,v)\mathbf{i}+y(u,v)\mathbf{j}+v\mathbf{k}
\]

where $\mathbf{i},\mathbf{j},\mathbf{k}$ is a right-handed orthonormal
basis of $E^{3}$ and $X_{u}\cdot X_{v}=0.$

We note that $\mathbf{i},\mathbf{j}$ is a basis of $E^{2}$ as a
subspace of $E^{3}.$

Definitions. $u,v$ are $\mathit{level}$ $\mathit{parameters}$ of
the surface.

$X_{0}$$=$$x(u,v)\mathbf{i}+y(u,v)\mathbf{j}$ is the $\mathit{principal}$
$\mathit{part}$ of the surface.

It follows that
\[
X(u,v)=X_{0}+v\mathbf{k}.
\]
 We need to define the function $X(u,v$) in order to find the function
giving the embedding of the given metric $\omega$. The embedding
function will be 
\[
\bar{X}(\bar{u},\bar{v})=X(u,v)
\]
 where the change of parameters from the barred parameters to the
unbarred parameters will be given by solving the pair of initial value
problems 
\begin{equation}
\sqrt{\bar{G}(\bar{u},\bar{v})}f_{\bar{u}}+\sqrt{\frac{f_{\bar{u}}^{2}}{1-f_{\bar{u}}^{2}}}f_{\bar{v}}=0\label{eq:IVP 1}
\end{equation}

with $C^{1}$ initial values $f(\bar{u},0)=\bar{h}(\bar{u})$ with
$\bar{h}'(\bar{u})$$\neq0$; 
\begin{equation}
\sqrt{\frac{f_{\bar{u}}^{2}}{1-f_{\bar{u}}^{2}}}g_{\bar{u}}\sqrt{\bar{G}(\bar{u},\bar{v})}-g_{\bar{v}}=0\label{IVP2}
\end{equation}

with $C^{1}$ initial values $g(\bar{u},0)=\bar{k}(\bar{u})$ with
$\bar{k}'(\bar{u})\neq0.$

Remark. We may choose the initial values to be not of class $C^{2}$
so that the solutions of the PDEs will not have power series in a
neighborhood of the initial point and hence will not be real analytic.
The isometric embedding will then be not real analytic. This will
be true even though the given metric is real analytic. An example
is 
\[
\omega=d\bar{u}^{2}+(\cos^{2}\bar{u})d\bar{v}^{2}.
\]
 An easy calculation using the Theorema Egregium of Gauss shows that
the curvature of the metric is equal to one in a neighborhood of $(0,0).$

Using the solutions of the initial value problems we define a change
of parameters by
\begin{equation}
u=f(\bar{u},\bar{v}),\;v=g(\bar{u},\bar{v}).\label{eq:change of params}
\end{equation}

Claim. The change of parameters $(3)$ is invertible. For a proof
see the Appendix of calculations.

We now define the composite vector function
\[
\bar{X}(\bar{u},\bar{v})=X(f(\bar{u},\bar{v}),g((\bar{u},\bar{v})).
\]

By the change of parameters $\bar{X}(\bar{u},\bar{v})$ will be well-defined
once $X(u,v)$ is defined. We can then write 
\[
\bar{X}(\bar{u},\bar{v})=X(u,v).
\]

Suppose for now that the functions $x(u,v),\,y(u,v)$ are known. This
will enable us to compute the components of the metric of the immersion
$\bar{X}(\bar{u},\bar{v})$ induced by the ambient space $E^{3}$
and show that they are equal respectively to $1,0,\bar{G}$, i.e.,
to the components of the given metric $\omega.$

\section{Definition of System S.}

The system S, to play an important role in our proof of embedding,
is
\[
1=1f_{\bar{u}}^{2}+Gg_{\bar{u}}^{2}.
\]
\[
0=1f_{\bar{u}}f_{\bar{v}}+Gg_{\bar{u}}g_{\bar{v}}.
\]
\[
\bar{G}=1f_{\bar{v}}^{2}+Gg_{\bar{v}}^{2}.
\]
 $\bar{G}=\bar{G}(\bar{u},\bar{v})$ comes from the given metric to
be embedded. S is algebraically linear in unknowns $1,G$ with coefficients
given by the partial derivatives which are known from the inital value
problems solved above.

Theorem 1. S has an algebraic solution for $1,G$ in terms of the
partial derivatives and $\bar{G}.$

Proof. See the appendix of calculations subsection 5.1.

We can now solve S for $1,G(u,v)$ by Cramer's rule.

$G(u,v)$ is now determined in terms of $\bar{G}(\bar{u},\bar{v}).$
In fact, from the third equation, 
\[
G(u,v)=\frac{\bar{G}-f_{\bar{v}}^{2}}{g_{\bar{v}}^{2}}.
\]

The system S says that $1,0,\bar{G}$ and $1,0,G$ are components
of the same metric. Therefore we have

Theorem 2. 
\[
\omega=d\bar{u}^{2}+\bar{G}(\bar{u},\bar{v})d\bar{v}^{2}=du^{2}+G(u,v)dv^{2}.
\]

QED theorem 2.

\section{The function X(u,v) is well-defined.}

Recall that
\[
X(u,v)=x(u,v)\mathbf{i}+y(u,v)\mathbf{j}+v\mathbf{k.}
\]

We first put the metric of the Euclidean plane in the geodesic form
\[
\omega_{0}=du^{2}+G_{0}(u,v)dv^{2}.
\]
 This can always be accomplished by choosing an arbitrary $C^{1}$
base curve as referenced above. If the base curve is a segment of
the x axis of a standard coordinate system of $E^{2},$ then the metric
is in the canonical form 
\[
\omega_{0}=dx^{2}+dy^{2}.
\]

The straight lines $x=$ constant are the geodesics. Let 
\[
X_{0}=x(u,v)\mathbf{i}+y(u,v)\mathbf{j}
\]
 represent $E^{2}$ in terms of geodesic parameters $u,v$. Thus $x(u,v),y(u,v)$
are determined once the base curve is given. The metric induced by
the ambient space $E^{2}$ on $X_{0}$ is given by 
\[
E_{0}=X_{0,u}=(x_{u}\mathbf{i}+y_{u}\mathbf{j})^{2}=x_{u}^{2}+y_{u}^{2}.
\]
 
\[
F_{0}=X_{0,u}\cdot X_{0,v}=(x_{u}\mathbf{i}+y_{u}\mathbf{j})\cdot(x_{v}\mathbf{i}+y_{v}\mathbf{j})=x_{u}x_{v}+y_{u}y_{v}=0.
\]
 (since geodesic parameters are orthogonal)
\[
G_{0}=X_{0,v}^{2}=(x_{v}\mathbf{i+y_{v}\mathbf{j})^{2}=x_{v}^{2}+y_{v}^{2},}
\]
 where $E_{0}=1,$ $F_{0}=0,\,G_{0}>0.$

Then the Euclidean metric $\omega_{0}$ of $E^{2}$ has components
$\,E_{0},F_{0},$$\,G_{0}$ which satisfy the system $S_{0}$ defined
as 
\[
1=x_{u}^{2}+y_{u}^{2},
\]
\[
0=x_{u}x_{v}+y_{u}y_{v},
\]
\[
G_{0}=x_{v}^{2}+y_{v}^{2}.
\]

$G_{0}$ is known from the geodesic form $\omega_{0}$ above.

The system $S_{0}$ expresses the relation between the components
$1,1$ of the Euclidean metric relative to the coordinates $x,y$
and its components relative to the parameters $u,v$ . Therefore the
Euclidean metric can be written 
\[
\omega_{0}=dx^{2}+dy^{2}=du^{2}+G_{0}(u,v)dv^{2}.
\]
 Clearly the geodesics of $E^{2}$ are straight line segments given
by $v=$ constant. Since $x(u,v)$, $y(u,v$) are known by the definition
of geodesic parameters, therefore

Theorem. $X(u,v)$ is well-defined as 
\[
X(u,v)=x(u,v)\mathbf{i}+y(u,v)\mathbf{j}+v\mathbf{k}.
\]

This is an immersion of $u,v$ space in $E^{3}$. We now show that
it is a surface: Calculate the metric induced by the ambient space
$E^{3}.$
\[
E(u,v)=X_{u}^{2}=x_{u}^{2}+y_{u}^{2},
\]
\[
F(u,v)=X_{u}\cdot X_{v}=x_{u}x_{v}+y_{u}y_{v}=0,
\]
\[
G(u,v)=x_{v}^{2}+y_{v}^{2}+1.
\]

Thus, $G(u,v)=G_{0}+1,$ $\,E(u,v)=E_{0}(u,v)=1.$ 
\[
EG-F^{2}\neq0.
\]

Therefore we have the

Theorem 3. $X(u,v)$ is a surface.

\section{Proof of Main Theorem.}

Main theorem. The composite vector function defined as
\[
\bar{X}(\bar{u},\bar{v})=X(f(\bar{u},\bar{v}),g(\bar{u},\bar{v}))
\]

is an isometric embedding of the given metric
\[
\omega=\bar{E}(\bar{u},\bar{v})d\bar{u}^{2}+\bar{F}d\bar{u}d\bar{v}+\bar{G}(\bar{u},\bar{v})d\bar{v}^{2}
\]

in $E^{3}.$

Proof. We compute the components of the metric of the immersion $\bar{X}(\bar{u},\bar{v})$
induced by the ambient space $E^{3}$ and show that they are equal
respectively to $1,0,\bar{G},$ i.e., to the components of the given
metric $\omega.$
\[
\bar{X}_{\bar{u}}^{2}=(X_{u}u_{\bar{u}}+X_{v}v_{\bar{u}})^{2}=X_{u}^{2}u_{\bar{u}}^{2}+2X_{u}\cdot X_{v}u_{\bar{u}}v_{\bar{u}}+X_{v}^{2}v_{\bar{u}}^{2}
\]
\[
=X_{u}^{2}u_{\bar{u}}^{2}+X_{v}^{2}v_{\bar{u}}^{2}
\]
\[
=E(u,v)u_{\bar{u}}^{2}+G(u,v)v_{\bar{u}}^{2}
\]
\[
=(x_{u}^{2}+y_{u}^{2})u_{\bar{u}}^{2}+(x_{v}^{2}+y_{v}^{2}+1)v_{\bar{u}}^{2}.
\]
 But $E=E_{0}=1,\:G=G_{0}+1$ implies 
\[
\bar{X}_{\bar{u}}^{2}=E_{0}(u,v)u_{\bar{u}}^{2}+(G_{0}(u,v)+1)v_{\bar{u}}^{2}
\]
 by the first equation of system S and $v=g(\bar{u},\bar{v}).$ Therefore
\[
\bar{X}_{\bar{u}}^{2}=1.
\]
 
\[
\bar{X}_{\bar{u}}\cdot\bar{X}_{\bar{v}}=(X_{u}u_{\bar{u}}+X_{v}v_{\bar{u}})\cdot(X_{u}u_{\bar{v}}+X_{v}v_{\bar{v}})
\]

\[
=X_{u}^{2}u_{\bar{u}}u_{\bar{v}}+X_{v}^{2}v_{\bar{u}}v_{\bar{v}}
\]
\[
=(x_{u}^{2}+y_{u}^{2})u_{\bar{u}}u_{\bar{v}}+(x_{v}^{2}+y_{v}^{2}+1)v_{\bar{u}}v_{\bar{v}}
\]
\[
=E_{0}u_{\bar{u}}u_{\bar{v}}+Gv_{\bar{u}}v_{\bar{v}}
\]
\[
=u_{\bar{u}}u_{\bar{v}}+Gg_{\bar{u}}g_{\bar{v}}
\]
 using the change of parameters (\ref{eq:change of params}). The
second equation of S may be written 
\[
0=1u_{\bar{u}}u_{\bar{v}}+Gg_{\bar{u}}g_{\bar{v}}.
\]
 Thus, 
\[
\bar{F}(\bar{u},\bar{v})=0.
\]
 
\[
\bar{X}_{\bar{v}}^{2}=(X_{u}u_{\bar{v}}+X_{v}v_{\bar{v}})^{2}
\]
\[
=(x_{u}^{2}+y_{u}^{2})u_{\bar{v}}^{2}+(x_{v}^{2}+y_{v}^{2}+1)v_{\bar{v}}^{2}
\]
\[
=Eu_{\bar{v}}^{2}+G(u,v)v_{\bar{v}}^{2}
\]
\[
=Ef_{\bar{v}}^{2}+G(u,v)g_{\bar{v}}^{2}
\]
\[
=f_{\bar{v}}^{2}+G(u,v)g_{\bar{v}}^{2}=\bar{G(}\bar{u},\bar{v})
\]
 by system S. QED main theorem.

\section{Appendix of Calculations.}

\subsection{Proof of Theorem 1. }

Theorem 1. S has an algebraic solution for $1,G$ in terms of the
partial derivatives and $\bar{G}.$

Proof. In order to prove theorem 1 we must show that S is algebraically
consistent. This will be true if and only if the rank of the augmented
matrix of the system equals the rank of its coefficient matrix \cite{key-12}.

Calculation of the ranks of the augmented matrix and the coefficient
matrix:

Lemma. The Jacobian J of the parameter change (3) is not zero in some
neighborhood of the initial point.

\[
J=\left|\begin{array}{cc}
f_{\bar{u}} & f_{\bar{v}}\\
g_{\bar{u}} & g_{\bar{v}}
\end{array}\right|.
\]
 At the initial point $f_{\bar{u}}(\bar{u},0)=\bar{h}'(\bar{u})\neq0.$
We may take $\bar{h}'(\bar{u})>0.$

Using initial value problem 1 (\ref{eq:IVP 1}) we calculate $f_{\bar{v}}(\bar{u},0):$
\[
\sqrt{\bar{G}(\bar{u},0)}\bar{h}'(\bar{u})+\sqrt{\frac{(\bar{h}'(\bar{u})^{2}}{1-(\bar{h}'(\bar{u})^{2}}}f_{\bar{v}}(\bar{u},0)=0.
\]

Cancelling $\bar{h}'(\bar{u})$ we obtain 
\[
\sqrt{\bar{G}(\bar{u},0)}+\sqrt{\frac{1}{1-(\bar{h}'(\bar{u})^{2}}}f_{\bar{v}}(\bar{u},0)=0.
\]
\[
f_{\bar{v}}(\bar{u},0)=-\sqrt{\bar{G}(\bar{u},0)}\sqrt{1-(\bar{h}'(\bar{u})^{2}}.
\]
 Next we calculate $g_{\bar{v}}$ using the solution of the second
initial value problem (\ref{IVP2}). At the initial point this becomes
\[
g_{\bar{v}}=\sqrt{\frac{(\bar{h}'(\bar{u})^{2}}{1-(\bar{h}'(\bar{u})^{2}}}g_{\bar{u}}\sqrt{\bar{G}(\bar{u},0)}=\sqrt{\frac{(\bar{h}'(\bar{u})^{2}}{1-(\bar{h}'(\bar{u})^{2}}}\bar{k}'(\bar{u})\sqrt{\bar{G}(\bar{u},0)}.
\]
 Therefore, 
\[
J(\bar{u},0)=\left|\begin{array}{cc}
f_{\bar{u}} & -\sqrt{\bar{G}(\bar{u},0)}\sqrt{1-(\bar{h}'(\bar{u})^{2}}\\
g_{\bar{u}} & \sqrt{\frac{(\bar{h}'(\bar{u})^{2}}{1-(\bar{h}'(\bar{u})^{2}}}\bar{k}'(\bar{u})\bar{G}(\bar{u},0)
\end{array}\right|
\]
\[
=\left|\begin{array}{cc}
\bar{h}'(\bar{u}) & -\sqrt{\bar{G}(\bar{u},0)}\sqrt{1-(\bar{h}'(\bar{u})^{2}}\\
\bar{k}'(\bar{u}) & \sqrt{\frac{(\bar{h}'(\bar{u})^{2}}{1-(\bar{h}'(\bar{u})^{2}}}\bar{k}'(\bar{u})\bar{G}(\bar{u},0)
\end{array}\right|
\]
\[
=\sqrt{1-(\bar{h}'(\bar{u}))^{2}}\left|\begin{array}{cc}
\bar{h}'(\bar{u}) & -\sqrt{\bar{G}(\bar{u},0)}\\
\bar{k}'(\bar{u}) & \bar{h}'(\bar{u})\frac{\bar{k}'(\bar{u})\bar{G}(\bar{u},0)}{1-(\bar{h}'(\bar{u}))^{2}}
\end{array}\right|.
\]

This can be made positive by taking $\bar{h}'(\bar{u})$ and $\bar{k}'(\bar{u})$
small positive. QED Lemma.

By the Lemma the parameter change is locally one-to-one, bicontinuous
and continuously differentiable, i.e. , is locally invertible.

We now compute the ranks of the cefficient and augmented matrices
of system S.

Claim. The rank of the augmented matrix is two.

Proof of Claim. 
\[
\left|\begin{array}{ccc}
f_{\bar{u}}^{2} & g_{\bar{u}}^{2} & 1\\
f_{\bar{u}}f_{\bar{v}} & g_{\bar{u}}g_{\bar{v}} & 0\\
f_{\bar{v}}^{2} & g_{\bar{v}}^{2} & \bar{G}(\bar{u},\bar{v})
\end{array}\right|
\]
\[
=\left|\begin{array}{cc}
f_{\bar{u}}f_{\bar{v}} & g_{\bar{u}}g_{\bar{v}}\\
f_{\bar{v}}^{2} & g_{\bar{v}}^{2}
\end{array}\right|+\bar{G}(\bar{u},\bar{v})\left|\begin{array}{cc}
f_{\bar{u}}^{2} & g_{\bar{u}}^{2}\\
f_{\bar{u}}f_{\bar{v}} & g_{\bar{u}}g_{\bar{v}}
\end{array}\right|
\]
\[
=f_{\bar{v}}g_{\bar{v}}\left|\begin{array}{cc}
f_{\bar{u}} & g_{\bar{u}}\\
f_{\bar{v}} & g_{\bar{v}}
\end{array}\right|+\bar{G}(\bar{u},\bar{v})f_{\bar{u}}g_{\bar{u}}\left|\begin{array}{cc}
f_{\bar{u}} & g_{\bar{u}}\\
f_{\bar{v}} & g_{\bar{v}}
\end{array}\right|
\]
\[
=J(f_{\bar{v}}g_{\bar{v}}+\bar{G}(\bar{u},\bar{v})f_{\bar{u}}g_{\bar{u}}.
\]
 Now substitute from the PDEs of the initial value problems to obtain
zero. Therefore the rank of the augmented matrix is at most two. That
the rank of the augmented matrix is at least two follows from the
Lemma and the fact that the determinant has two by two minors that
are not zero. Therefore the rank of the augmented matrix is two. QED
Claim.

Claim'. The rank of the coefficient matrix is two.

Proof of Claim'. The coefficient matrix is
\[
\left[\begin{array}{cc}
f_{\bar{u}}^{2} & g_{\bar{u}}^{2}\\
f_{\bar{u}}f_{\bar{v}} & g_{\bar{u}}g_{\bar{v}}\\
f_{\bar{v}}^{2} & g_{\bar{v}}^{2}
\end{array}\right].
\]
\[
\left|\begin{array}{cc}
f_{\bar{u}}^{2} & g_{\bar{u}}^{2}\\
f_{\bar{u}}f_{\bar{v}} & g_{\bar{u}}g_{\bar{v}}
\end{array}\right|=f_{\bar{u}}g_{\bar{u}}\left|\begin{array}{cc}
f_{\bar{u}} & g_{\bar{u}}\\
f_{\bar{v}} & g_{\bar{v}}
\end{array}\right|.
\]
 At the initial point, where $\bar{v}=0$, this becomes 
\[
J(\bar{u},0)h'(\bar{u})k'(\bar{u})\neq0.
\]
 Therefore in a neighborhood of the initial point the rank of the
coefficient matrix is two. QED Claim'. QED Theorem 1.

\subsection{Proof of the existence of the inverse of the prameter change.}

Lemma. The Jacobian of the parameter change is not zero locally.

Proof. Define 
\begin{align*}
\lambda(\bar{u},\bar{v}) & =\sqrt{\frac{f_{\bar{u}}^{2}}{1-f_{\bar{u}}^{2}}}.
\end{align*}

Then the PDEs can be written 
\[
\sqrt{\bar{G}(\bar{u},\bar{v})}f_{\bar{u}}+\lambda(\bar{u},\bar{v})f_{\bar{v}}=0.
\]
\[
\lambda(\bar{u},\bar{v})g_{\bar{u}}\sqrt{\bar{G}(\bar{u},\bar{v})}-g_{\bar{v}}=0.
\]
 
\[
J=\left|\begin{array}{cc}
f_{\bar{u}} & f_{\bar{v}}\\
g_{\bar{u}} & g_{\bar{v}}
\end{array}\right|=\left|\begin{array}{cc}
f_{\bar{u}} & -\frac{\sqrt{\bar{G}(\bar{u},\bar{v})}f_{\bar{u}}}{\lambda}\\
g_{\bar{u}} & \sqrt{\bar{G}(\bar{u},\bar{v})\lambda}g_{\bar{u}}
\end{array}\right|=f_{\bar{u}}g_{\bar{u}}\sqrt{\bar{G}}\left|\begin{array}{cc}
1 & -\lambda^{-1}\\
1 & \lambda
\end{array}\right|.
\]

By a suitable choice of the initial conditions we can make $J\neq0$
in a neithborhood of the initial point. Choose 
\[
\bar{h}'(\bar{u}')>0,
\]
\[
\bar{k}'(\bar{u})>0.
\]
 Thus $J>0$ taking care to choose $\bar{h}'(\bar{u}')$ small enough
that $\lambda$ is real. By a text book result the mapping is locally
invertible. QED.

\section{Topic for further investigation. Informal Remarks.}

The existence of isometric embeddings of $C^{1}$ metrics which are
$C^{2}.$ Compatibility conditions are needed for a proof since the
systems are overdetermined.

The principal part of the embedded surface, defined in the introduction,
is two-dimensional and serves to determine a surface in a three-dimensional
space. This essentially reduces a three-dimensional problem to a two-dimensional
problem. It relates to an interesting discussion of Brian Greene in
a chapter of his book (The Fabric of the Cosmos) entitled 'Is the
Universe a Hologram?' ``Whereas Plato envisioned common perceptions
as revealing a mere shadow of reality, the holographic principle concurs,
but turns the metaphor on its head. The shadows$-$the things that
are flattened out and hence live on a lower dimensional surface $-$
are real, while what seems to be the more richly structured, higher-dimensional
entities (us; the world around us) are evanescent projections of the
shadows.

{[}4{]} Brian Greene. The Fabric of the Cosmos. Alfred A. Knopf, New
York, 2004.
\end{document}